\documentclass[a4paper,fleqn]{cas-sc}

\usepackage[authoryear]{natbib}

\usepackage{amsmath,bm}
\usepackage{graphicx}
\usepackage{algorithm}
\usepackage{algpseudocode}
\usepackage{hyperref}

\newcommand\ktxt{\textcolor{black}}

\newcommand{\orcidicon}[1]{%
	\href{https://orcid.org/#1}{%
		\includegraphics[height=1.5ex]{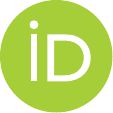}%
	}%
}

\def\tsc#1{\csdef{#1}{\textsc{\lowercase{#1}}\xspace}}
\tsc{WGM}
\tsc{QE}

\begin{document}
\let\WriteBookmarks\relax
\def\floatpagepagefraction{1}
\def\textpagefraction{.001}

\shorttitle{}    

\shortauthors{}  

\title [mode = title]{Mode-realigned pointwise interpolation (MRPWI) for efficient POD-Galerkin parametric reduced-order models}  

\author[1]{\ktxt{Lei Du} \orcidicon{0009-0000-2079-7019}}
\credit{Conceptualization, Data curation, Formal analysis, Investigation, Methodology, Software, Validation, Visualization, Writing – original draft, Writing – review \& editing}

\author[1,2]{\ktxt{Shengqi Zhang} \orcidicon{0000-0001-8273-7484}}
\cormark[1]
\cortext[1]{Corresponding author}
\ead{szhang@eitech.edu.cn}
\credit{Conceptualization, Formal analysis, Funding acquisition, Investigation, Methodology, Project administration, Resources, Supervision, Writing – review \& editing}

\affiliation[1]{organization={Zhejiang Key Laboratory of Industrial Intelligence and Digital Twin, Eastern Institute of Technology},
	city={Ningbo},
	postcode={315200}, 
	country={China}}
\affiliation[2]{organization={Ningbo Institute of Digital Twin, Eastern Institute of Technology},
	city={Ningbo},
	postcode={315200}, 
	country={China}}
	
\begin{abstract}
As a cornerstone of reduced-order modeling, the POD-Galerkin framework has garnered widespread attention and remains one of the most widely adopted approaches. Constructing POD-Galerkin PROMs involves integrating this framework with advanced interpolation techniques to obtain POD modes at target (unseen) parameters. While Grassmann manifold interpolation (GMI) serves as an accurate baseline, mode-realigned pointwise interpolation (MRPWI) is proposed to develop highly efficient PROMs that maintain comparable accuracy. Notably, the MRPWI employs a two-step mode realignment procedure, consisting of sign alignment and rotation alignment, to effectively synchronize the POD modes. Demonstration and evaluation of the constructed POD-Galerkin PROMs are conducted by examining flow over a cylinder. These models exhibit high fidelity in comparison to direct numerical simulation and standard POD-Galerkin ROMs. PROMs constructed via MRPWI achieve accuracy comparable to those using GMI, while providing significantly higher computational efficiency.
\end{abstract}



\maketitle

\section{Introduction}\label{intro}
Parametric reduced-order models (PROMs), reduced-order models (ROMs) tailored for parameter-dependent problems, have emerged as a promising tool to replace full-order models (FOMs) in science and engineering \citep{brunton2020machine,mendez2023data}. A prominent class of such problems arises from parametric partial differential equations, which govern a wide range of physical phenomena, yet whose numerical simulation often leads to high-dimensional dynamical systems that are computationally prohibitive to solve repeatedly under varying parameters \citep{hinze2023model}. This motivates the development of PROMs that reduce the computational burden of FOMs while maintaining accuracy and efficiency.

Proper orthogonal decomposition (POD) is one of the most widely employed methods for modal decomposition and reduced-order modeling, serving as a foundation and inspiration for a broad class of related modal decomposition techniques \citep{taira2017modal}. By identifying a low-dimensional subspace that optimally approximates high-dimensional data, POD yields compact, low-rank representations that accurately capture the spatiotemporal evolution of complex dynamical systems \citep{brunton2022data}. \citet{lumley1970stochastic} first introduced POD to the fluid mechanics community for the analysis of turbulent flows. In practice, the widely adopted method of snapshots proposed by \citet{sirovich1987turbulence} determines the POD modes via eigendecomposition of the snapshot correlation matrix, which is computationally equivalent to performing a singular value decomposition (SVD) of the data matrix. In either formulation, POD extracts the most energetically dominant features from the snapshot ensemble. Its generality and interpretability have enabled broad application across disciplines, including fluid dynamics \citep{holmes2012turbulence}, acoustics \citep{feeny2002proper}, and image processing \citep{sirovich1987low}. Additionally, several variants have been developed for specific purposes, including balanced POD \citep{willcox2002balanced,rowley2005model}, spectral POD \citep{towne2018spectral}, and multi-scale POD \citep{mendez2019multi}.

POD has been integrated with various methods to generate PROMs, including interpolation and regression approaches. A global POD basis, constructed from snapshots spanning the entire parameter space, may fail to optimally represent the features at any specific parameter, limiting accuracy in problems with strong parameter dependence \citep{benner2021system}. Therefore, methods that seek locally optimal POD modes at target parameters are preferred over global approaches \citep{benner2015survey}. Interpolation-based methods, such as subspace angle interpolation (SAI) \citep{lieu2004parameter} and manifold interpolation \citep{amsallem2008interpolation, benner2021system}, are among the most established, enabling the systematic development of POD PROMs across the parameter space. Regression-based approaches, such as Gaussian process subspace prediction \citep{zhang2022gaussian} and neural networks employing multilayer perceptrons \citep{hess2020comparison}, kernel-based shallow neural networks \citep{kapadia2024active}, and factorized Fourier neural operators \citep{fanaskov2025deep}, provide alternatives to interpolation for constructing locally optimal subspaces.

Data-driven PROMs offer advantages in implementation simplicity and an equation-free nature, but suffer from limited generalizability and lack of physical consistency \citep{pan2025physics,prakash2025nonintrusive}. Although hybrid approaches that combine physics-based and data-driven PROMs have been explored to address this limitation \citep{chen2021physics,pan2025physics}, they still rely on a global POD basis and remain in early stages of development. Projection-based methods, on the other hand, are grounded in the underlying governing equations, making them typically more systematic and theoretically rigorous \citep{benner2021system, mendez2023data}. While established interpolation methods such as SAI and Grassmann manifold interpoaltion (GMI) have demonstrated satisfactory accuracy in constructing POD-Galerkin PROMs, the development of approaches capable of supporting high-order interpolation with greater computational efficiency remains an open challenge. To this end, mode-realigned pointwise interpolation (MRPWI) is proposed to recover POD modes at target parameters, enabling POD-Galerkin PROMs that retain the accuracy benefits of high-order interpolation while achieving enhanced efficiency. The method follows a two-step mode realignment procedure, comprising sign alignment and rotation alignment derived from Kasner's pseudo-angle. The resulting PROMs are evaluated on flow over a cylinder, with DNS, standard POD-Galerkin ROMs and those constructed with GMI serving as references.

The paper proceeds as follows. POD-Galerkin ROMs, covering both POD and Galerkin projection, are presented in Section~\ref{sec:pod_roms}. Section~\ref{sec:pod_proms} presents the interpolation methods, namely GMI and MRPWI, and their role in constructing POD-Galerkin PROMs. A test case involving flow over a cylinder is then examined in Section~\ref{sec:test_case}, followed by conclusions in Section~\ref{sec:conclusions}.

\section{POD-Galerkin reduced-order models}\label{sec:pod_roms}
\subsection{Proper orthogonal decomposition}\label{sec:pod}
POD is a widely adopted dimensionality reduction technique that extracts dominant coherent structures from a set of snapshots, yielding energy-optimal spatial modes \citep{kutz2016dynamic,taira2020modal,callaham2023multiscale}.

Initially, the data matrix $Q$ is formed from the snapshots,
\begin{equation}
	Q = \left[\bm{q}_{1}, \bm{q}_{2}, \cdots, \bm{q}_{N_s}\right] \in \mathbb{R}^{N_n \times N_s}
	\label{eq:data_matrix}
\end{equation}
where each snapshot $\bm{q}_{j}$ represents the state of the high-dimensional dynamical system at time $t_j$ and is reshaped into a column vector. Incompressible flows are considered, for which
\begin{equation}
	\bm{q}_{j} = \begin{bmatrix} \bm{u}(\bm{x},t_j) \\ \bm{p}(\bm{x},t_j) \end{bmatrix} \in \mathbb{R}^{N_n}.
\end{equation}
Here, $\bm{q}_j$ comprises the flow variables, namely velocity and pressure, $\bm{x}$ denotes the spatial coordinate, and $j$ indexes the snapshots. In such problems, the data matrix $Q$ is typically tall and skinny, as the spatial dimension $N_n$ is significantly larger than the number of snapshots $N_s$. 

To accommodate data defined on a non-uniform mesh, the weighted inner product is introduced as
\begin{equation}
	\langle\bm{q}_1, \bm{q}_2\rangle = \bm{q}_1^T \bm{W} \bm{q}_2,
	\label{eq:inner_product}
\end{equation}
where superscript $T$ denotes the matrix transpose, and $\bm{W}$ is a weight matrix whose entries reflect the numerical quadrature weights associated with the local cell volumes of the spatial discretization.

Once weighted, the data matrix admits a rank-$N_r$ truncated SVD,
\begin{equation}
	\bm{W}^{1/2} \bm{Q} \approx \bm{U} \bm{\Sigma} \bm{V}^T,
\end{equation}
where the non-negative singular values are given by the diagonal entries of $\bm{\Sigma} \in \mathbb{R}^{N_r \times N_r}$, and $\bm{U} \in \mathbb{R}^{N_n \times N_r}$ and $\bm{V} \in \mathbb{R}^{N_s \times N_r}$ contain the left and right singular vectors, which encode the spatial and temporal correlations inherent in the data, respectively. Consequently, the POD modes are obtained as,
\begin{equation}
	\bm{\Phi} = \bm{W}^{-1/2} \bm{U} = [ \bm{\phi}^1, \bm{\phi}^2, \ldots, \bm{\phi}^{N_r} ] \in \mathbb{R}^{N_n \times N_r},
\end{equation}
which satisfy the orthonormality condition $\bm{\Phi}^T \bm{W} \bm{\Phi} = \bm{I}$ with respect to the weighted inner product defined in Eq.~(\ref{eq:inner_product}) \citep{massaro2023flow}.

\subsection{Galerkin projection}\label{sec:pod-galerkin}
The POD-Galerkin framework combines the dimensionality reduction method, POD, with Galerkin projection to construct a low-dimensional dynamical system, i.e., the POD-Galerkin ROM, capable of accurately modeling full-order dynamical systems.

In the context of the incompressible Navier-Stokes equations, explicitly including the pressure term within the formulation is known to significantly enhance the accuracy and stability of the resulting ROMs, as demonstrated in various studies \citep{noack2005need, bergmann2009enablers, holmes2012turbulence}. Accordingly, the velocity $\bm{u}$ and pressure $p$ are approximated as a truncated linear superposition of orthonormal spatial modes $\bm{\phi}_l$ with temporal coefficients $\alpha_l$,
\begin{equation}
	\begin{bmatrix} \bm{u}(\bm{x},t) \\ \bm{p}(\bm{x},t) \end{bmatrix} \approx \sum_{l=1}^{N_r} \alpha_l(t) \begin{bmatrix} \bm{\varphi}_l(\bm{x}) \\ \bm{\psi}_l(\bm{x}) \end{bmatrix} = \sum_{l=1}^{N_r} \alpha_l(t) \bm{\phi}_l(\bm{x}),
	\label{eq:approximation}
\end{equation}
where $\bm{\phi}_l$ denotes the $l$-th spatial mode, composed of the velocity component $\bm{\varphi}_l$ and the pressure component 
$\bm{\psi}_l$ \citep{lorenzi2016pod}. These modes are obtained via POD applied to the data matrix in Eq.~(\ref{eq:data_matrix}), as discussed in detail in Section~\ref{sec:pod}.

A Galerkin projection of the Navier-Stokes momentum equation onto the finite-dimensional subspace spanned by the spatial modes yields \citep{taira2020modal},
\begin{equation}
	\langle \bm{\varphi}_l,\bm{u}_t \rangle = \langle \bm{\varphi}_l, -\bm{u} \cdot \nabla \bm{u} \rangle + \langle \bm{\varphi}_l, \nu\nabla^2\bm{u} - \nabla p \rangle.
	\label{eq:projection}
\end{equation}
Substituting the expansion from Eq.~(\ref{eq:approximation}) into Eq.~(\ref{eq:projection}) leads to a system of $N_r$ coupled ordinary differential equations that govern the evolution of the temporal coefficients,
\begin{equation}
	\sum_{m=1}^{N_r} M_{lm} \dot{\alpha}_m = \sum_{m=1}^{N_r} \sum_{n=1}^{N_r} Q_{lmn} \alpha_{m} \alpha_{n} + \sum_{m=1}^{N_r} L_{lm} \alpha_{m},
	\label{eq:galerkin}
\end{equation}
which constitutes the final POD-Galerkin ROM. Here, the system operators are defined as $M_{lm} = \langle \bm{\varphi}_{l},\bm{\varphi}_m \rangle$, $Q_{lmn} = \langle \bm{\varphi}_{l}, - \bm{\varphi}_m \cdot \nabla \bm{\varphi}_n \rangle$, and $L_{lm} = \langle \bm{\varphi}_{l}, \nu \nabla^2 \bm{\varphi}_m - \nabla \bm{\psi}_m \rangle$, representing the mass, quadratic (convective), and linear (viscous and pressure gradient) terms, respectively.

\section{POD-Galerkin parameteric reduced-order modes}\label{sec:pod_proms}
The key to constructing POD-Galerkin PROMs is to determine the POD modes $\overline{\mathit{\Phi}}$ at the target parameter $\overline{\gamma}$. To this end, interpolation methods, namely GMI and MRPWI, are employed to approximate $\overline{\mathit{\Phi}}$ from the neighboring POD modes. In addition, the selection criteria for both the $N_p$ interpolation cases and the reference case are identical to those in \citet{du2026interpolation}. Algorithm~\ref{alg} details the steps for constructing POD-Galerkin PROMs.
\begin{algorithm}[H]
	\caption{Procedure for constructing POD-Galerkin PROMs}
	\label{alg}
	\begin{algorithmic}[1]
		\Require Flow data $\{Q_j \in \mathbb{R}^{N_n \times N_s}\}_{j=1}^{N_c}$ associated with system parameters $\{\gamma_j\}_{j=1}^{N_c}$.
		\Ensure Predicted state $\bm{q}$ at the target parameter $\overline{\gamma}$.
		\State Extract the POD modes $\{\mathit{\Phi}_j \in \mathbb{R}^{N_n \times N_r}\}_{j=1}^{N_c}$ for the available cases with parameters $\{\gamma_j\}_{j=1}^{N_c}$. Section~\ref{sec:pod}
		\State Determine $N_p$ neighboring cases from the $N_c$ available cases and select one as the reference ($j=0$).
		\State Approximate the POD modes $\overline{\mathit{\Phi}}$ at $\overline{\gamma}$ via interpolation methods such as GMI and MRPWI from the modes $\{\mathit{\Phi_j}\}_{j=1}^{N_p}$. Sections~\ref{sec:gmi}, \ref{sec:mrpwi}
		\State Predict the system evolution according to Eq.~(\ref{eq:galerkin}). Section~\ref{sec:pod-galerkin}
	\end{algorithmic}
\end{algorithm}

\subsection{Grassmann manifold interpolation}\label{sec:gmi}
GMI enables subspace interpolation while preserving the underlying geometric structure, making it well suited for the construction of PROMs \citep{benner2021system}. As the method has a well-established theoretical foundation and geometric interpretation \citep{amsallem2008interpolation,amsallem2010interpolation,bendokat2024grassmann,du2025interpolation,du2026interpolation}, only the resulting procedure is presented here for brevity.

POD modes $\{\mathit{\Phi}_j\}_{j=1}^{N_p}$, associated with the system parameters $\{\gamma_j\}_{j=1}^{N_p}$, are considered as subspaces that lie on the Grassmann manifold. Accordingly, GMI is employed to determine the POD modes $\overline{\mathit{\Phi}}$ at the target parameter $\overline{\gamma}$, which also lie on the manifold. At the first step, a group of reference modes $\mathit{\Phi}_0$ is selected from $\{\mathit{\Phi}_j\}_{j=1}^{N_p}$. Following this, the modes $\{\mathit{\Phi}_j\}_{j=1}^{N_p}$ are mapped to the tangent space of the Grassmann manifold at the reference modes $\mathit{\Phi}_0$ via the logarithmic mapping,
\begin{subequations}
	\begin{align}
		(I-\mathit{\Phi}_{0} \mathit{\Phi}_{0}^{*}) \mathit{\Phi}_{j} (\mathit{\Phi}_{0}^{T} \mathit{\Phi}_{j})^{-1} &= M_{j} \mathit{\Xi}_j N_{j}^{T} \quad (\text{thin SVD}), \\
		\mathit{\Gamma}_{j} &= M_{j} \text{tan}^{-1} \mathit{\Xi}_{j} N_{j}^{T}.
	\end{align}
\end{subequations}
Here, $\mathit{\Gamma}_{j}$ denotes the projection onto the tangent space. The corresponding $\overline{\mathit{\Gamma}}$ at the target parameter $\overline{\gamma}$ is obtained through standard Euclidean interpolation of the set $\{\mathit{\Gamma}_{j}\}_{j=1}^{N_p}$,
\begin{equation}
	\overline{\mathit{\Gamma}} = \sum\limits_{j=1}\limits^{N_p} \sigma_{j} \mathit{\Gamma}_{j},
	\label{eq:gmi3}
\end{equation}
where $\sigma_j$ are the Lagrange interpolation weights. Finally, the interpolated $\overline{\mathit{\Gamma}}$ are mapped back onto the Grassmann manifold via the exponential mapping, 
\begin{subequations}
	\begin{align}
		\overline{\mathit{\Gamma}} & = \overline{M} \ \overline{\mathit{\Xi}} \ \overline{N}^{T} \quad (\text{thin SVD}), 
		\\
		\overline{\mathit{\Phi}} & = \left( \mathit{\Phi}_{0} \overline{N} \text{cos} \overline{\mathit{\Xi}} + \overline{M} \text{sin} \overline{\mathit{\Xi}} \right) \overline{N}^{T},
	\end{align}
\end{subequations}
thus yielding the POD modes $\overline{\mathit{\Phi}}$ at the target parameter $\overline{\gamma}$.

\subsection{Mode-realigned pointwise interpolation}\label{sec:mrpwi}
MRPWI is a highly accurate, efficient, and simple method for mode interpolation in the construction of PROMs. It aligns all modes associated with different system parameters to a set of reference modes selected from them, and subsequently performs pointwise interpolation on the aligned modes to obtain those at the target parameter.  Although it has been proposed and demonstrated in \citet{du2025interpolation, du2026interpolation}, it cannot be directly applied to POD modes. Therefore, MRPWI is modified to adapt to the POD-Galerkin framework for constructing POD-Galerkin PROMs.

Once the POD modes $\{\mathit{\Phi}_j\}_{j=1}^{N_p}$ corresponding to the system parameters $\{\gamma_j\}_{j=1}^{N_p}$ have been determined, a set of modes is chosen from $\{\mathit{\Phi}_j\}_{j=1}^{N_p}$ to serve as the reference modes, denoted by $\mathit{\Phi}_0$. It is noted that $\mathit{\Phi}_j = [ \bm{\phi}_j^1, \bm{\phi}_j^2, \ldots, \bm{\phi}_j^{N_r} ] \in \mathbb{R}^{N_n \times N_r}$ consists of POD modes as its column vectors, where $\bm{\phi}_j^k \in \mathbb{R}^{N_n}$ for $k = 1,\ldots, N_r$. To begin the procedure, each set of POD modes $\mathit{\Phi}_j$ is arranged into complex vectors as follows,
\begin{equation}
	\hat{\bm{\phi}}_j^{\,k} =
	\begin{cases}
		\bm{\phi}_j^{1} + i \bm{0}, \\ 
		\bm{\phi}_j^{2k-2} + i\,\bm{\phi}_j^{2k-1}, \quad k = 2, \ldots, (N_r+1)/2,
	\end{cases}
\end{equation}
where $i$ denotes the imaginary unit and $\hat{\mathit{\Phi}}_j \in \mathbb{C}^{N_n \times (N_r+1)/2}$, with $N_r$ assumed to be odd. This form is well suited for mode realignment using Kasner's pseudo-angle, facilitating consistent interpolation across different parameters. Following this, the modes $\{\mathit{\hat{\Phi}}_j\}_{j=1}^{N_p}$ are realigned to the reference modes $\hat{\mathit{\Phi}}_0$ columnwise through a two-step procedure. The first step is the sign alignment, given by,
\begin{equation}
	\bm{\hat{\phi}}_j^{k} 
	= \text{sgn}\left(\text{Re}\left[\bm{\hat{\phi}}_0^{k}\right]^T \text{Re}\left[\bm{\hat{\phi}}_j^{k}\right]\right) \text{Re}\left[\bm{\hat{\phi}}_j^{k}\right] 
	+ i \cdot \text{sgn}\left(\text{Im}\left[\bm{\hat{\phi}}_0^{k}\right]^T \text{Im}\left[\bm{\hat{\phi}}_j^{k}\right]\right)\text{Im}\left[\bm{\hat{\phi}}_j^{k}\right].
\end{equation}
Here, $\mathrm{sgn}(\cdot)$ is the sign function, while $\mathrm{Re}[\cdot]$ and $\mathrm{Im}[\cdot]$ represent the real and imaginary parts, respectively. Next, the modes are rotated to align with the reference modes,
\begin{equation}
	\bm{\hat{\phi}}_j^{k} \mapsto \bm{\hat{\phi}}_j^{k} e^{-i\varphi_j^k}, \quad j=1,\ldots,N_p, \quad k=1,\ldots,(N_r+1)/2,
\end{equation}
where $\varphi_j^k = \varphi_j^k(\bm{\hat{\phi}}_0^k, \bm{\hat{\phi}}_j^k)$ is the Kasner's pseudo-angle, defined as,
\begin{equation}
	e^{i\varphi_j^k(\bm{\hat{\phi}}_0^k, \bm{\hat{\phi}}_j^k)} = \frac{\langle \bm{\hat{\phi}}_0^k, \bm{\hat{\phi}}_j^k \rangle}{|\langle \bm{\hat{\phi}}_0^k, \bm{\hat{\phi}}_j^k \rangle|}
	\quad \Rightarrow \quad
	\varphi_j^k(\bm{\hat{\phi}}_0^k, \bm{\hat{\phi}}_j^k) = -i \ln \frac{\langle \bm{\hat{\phi}}_0^k, \bm{\hat{\phi}}_j^k \rangle}{|\langle \bm{\hat{\phi}}_0^k, \bm{\hat{\phi}}_j^k \rangle|}.
\end{equation}
The angle $\varphi_j^k$ quantifies the relationship between the complex vectors $\bm{\hat{\phi}}_0^k$ and $\bm{\hat{\phi}}_j^k$ and lies in the range $\varphi_j^k \in [-\pi, \pi]$, where the notation $\langle \bm{\hat{\phi}}_0^k, \bm{\hat{\phi}}_j^k \rangle$ represents the complex inner product $(\bm{\hat{\phi}}_0^k)^* \bm{\hat{\phi}}_j^k$, with $*$ denoting the conjugate transpose. These aligned modes $\{\mathit{\hat{\Phi}}_j\}_{j=1}^{N_p}$ are then utilized to obtain the counterparts $\overline{\mathit{\hat{\Phi}}}$ at the target parameter $\overline{\gamma}$ via Lagrange interpolation,
\begin{equation}
	\overline{\mathit{\hat{\Phi}}} = \sum\limits_{j=1}\limits^{N_p} \sigma_{j} \mathit{\hat{\Phi}}_{j},
\end{equation}
where $\sigma_j$ are the interpolation weights, following Eq.~(\ref{eq:gmi3}). Finally, the POD modes $\overline{\mathit{\Phi}}$ at the target parameter $\overline{\gamma}$ are recovered as,
\begin{subequations}
	\begin{align}
		& \overline{\bm{\phi}}^{1} := \text{Re}\left[\bm{\overline{\hat{\phi}}}^{1}\right],  \\ 
		& \overline{\bm{\phi}}^{2k-2} := \text{Re}\left[\bm{\overline{\hat{\phi}}}^{k}\right], \quad k=2,3,\ldots,(N_r+1)/2, \\
		& \overline{\bm{\phi}}^{2k-1} := \text{Im}\left[\bm{\overline{\hat{\phi}}}^{k}\right], \quad k=2,3,\ldots,(N_r+1)/2.
	\end{align}
\end{subequations}

\section{Illustrative application: flow over a cylinder}\label{sec:test_case}
A canonical case, flow over a cylinder, is selected to demonstrate the POD-Galerkin PROMs, which are constructed with GMI and MRPWI. These PROMs are examined with respect to three key factors: the number of modes $N_r$, the system parameter interval $\mathrm{\Delta} Re$, and the number of neighbors $N_p$, with the target parameter set to $\overline{Re}=130$. 

\subsection{Flow configuration, governing equations, and numerical framework}
\begin{figure}
	\centering
	\includegraphics[width=.5\textwidth]{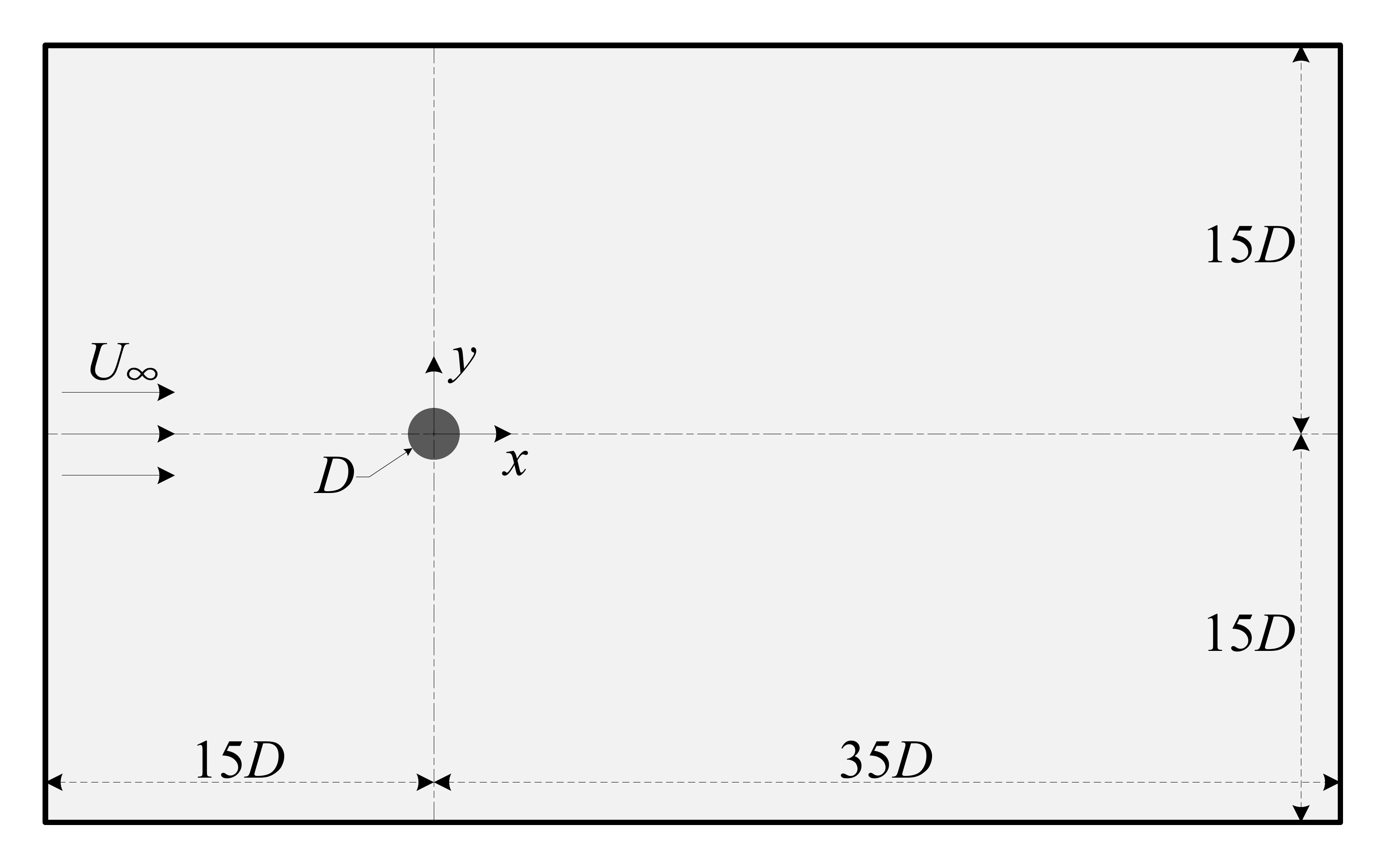}
	\caption{Flow over a cylinder: computational domain with cylinder diameter doubled for visualization.}
	\label{fig:domain}
\end{figure}

A schematic of the flow over a cylinder is presented in Fig.~\ref{fig:domain}, where the computational domain follows the layout described in \citet{du2025interpolation,du2026interpolation}. A Cartesian frame $\bm{x} = (x,y)^T$ is introduced with its origin at the cylinder center, where $x$ and $y$ indicate the streamwise and transverse coordinates, respectively. The unsteady, two-dimensional incompressible flow is described by the dimensionless Navier-Stokes system, consisting of the continuity and momentum equations,
\begin{subequations}
	\begin{align}
		\nabla \cdot \bm{u} & = 0, \\
		\partial \bm{u}/\partial t + \bm{u} \cdot \nabla \bm{u}  & = -\nabla p + Re^{-1} \ \nabla^2 \bm{u}.
	\end{align}
\end{subequations}
Here, $\bm{u}(\bm{x},t) = (u,v)^T$ and $p(\bm{x},t)$ denote the velocity and pressure, respectively, defined over $\bm{x} = (x,y)^T$ and time $t$. The Reynolds number is defined as $Re = U_\infty D / \nu$, where $U_\infty$, $D$, and $\nu$ denote the freestream velocity, characteristic length (cylinder diameter), and kinematic viscosity, respectively.

The computational domain spans $x = -15D$ to $x = 35D$ in the streamwise direction and $y = -15D$ to $y = 15D$ in the transverse direction. The initial condition is a spatially uniform velocity, i.e., $\bm{u} = (U_{\infty}, 0)$, from which unsteady vortex shedding emerges during the temporal evolution of the simulation. This same velocity profile is prescribed at the inlet ($x=-15D$), corresponding to a Dirichlet boundary condition. A natural boundary condition is enforced at the outlet ($x = +35D$), given by $\left[-p \bm{I} + \nu \nabla \bm{u} \right] \cdot \bm{n} = 0$, where $\bm{n}$ is the unit normal vector. For the top ($y = +15D$) and bottom ($y = -15D$) boundaries, periodic boundary conditions are applied. On the cylinder surface, a no-slip wall condition (i.e., a Dirichlet boundary condition with $\bm{u}=\bm{0}$) is imposed \citep{massaro2024global}.

The governing equations, namely the Navier-Stokes equations, are time-integrated via direct numerical simulation (DNS) with the spectral element code Nek5000 \citep{patera1984spectral,deville2002high}. The computational domain is discretized into 1472 spectral elements, each containing $8^2 = 64$ Gauss-Lobatto-Legendre (GLL) quadrature points. Consequently, in each coordinate direction within an element, 
the velocity and pressure fields are represented by seventh-order Lagrangian polynomials defined on the GLL quadrature points, following the $\mathbb{P}_N-\mathbb{P}_{N}$ formulation. Time advancement is carried out with a third-order backward difference/extrapolation scheme (BDF3/EXT3), in which viscous terms are treated implicitly via BDF3 and convective terms advanced explicitly through EXT3 \citep{callaham2023multiscale}, with a fixed time step of $\Delta \tau = 0.001$ throughout all simulations. Following these settings, seventeen Reynolds number cases, $Re \in \{$70, 80, 85, 100, 110, 115, 120, 125, 130, 135, 140, 145, 150, 160, 175, 180, 190$\}$, are simulated to evaluate the POD-Galerkin PROMs. Flow data generated from these simulations are recorded during the periodic vortex shedding regime at intervals of $\Delta t = 100\,\Delta \tau = 0.1$.

\subsection{Evaluating POD-Galerkin parameteric reduced-order modes}
At the target parameter, the POD-Galerkin PROMs are characterized by the relative $L^2$ error (RLE),
\begin{equation}
	\bar{\varepsilon} = \sqrt{\frac{\sum_{j=1}^{N_\tau} \displaystyle \left\| \bm{q}_j^{\text{ROM}/\text{PROM}} - \bm{q}_j^{\text{DNS}} \right\|_2^2}{\sum_{j=1}^{N_\tau} \displaystyle \left\| \bm{q}_j^{\text{DNS}} \right\|_2^2}}.
\end{equation}
In this definition, $\bm{q}_j$ corresponds to the $j$-th snapshot from DNS, ROM, or PROM, where DNS provides the reference, and the averaging is performed over $N_\tau$ snapshots.

\begin{figure}
	\centering
	\begin{minipage}{0.8\textwidth}
		\centering
		\includegraphics[width=\textwidth]{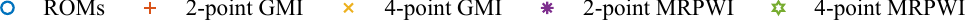}
	\end{minipage}  
	
	\vskip 0.3cm
	
	\begin{minipage}{0.8\textwidth}
		\centering
		\includegraphics[width=\textwidth]{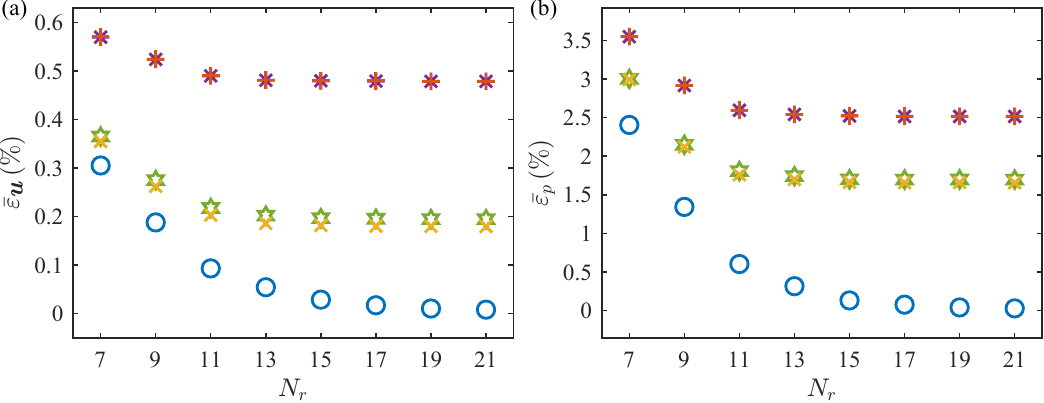}
	\end{minipage} 
	
	\caption{RLEs with respect to number of modes $N_r$ for (a) velocity and (b) pressure. (2-point: $Re \in \{120, 140\}$; 4-point: $Re \in \{100, 120, 140, 160\}$).}
	\label{fig:r}
\end{figure}

\begin{figure}
	\centering
	\begin{minipage}{0.8\textwidth}
		\centering
		\includegraphics[width=\textwidth]{pod_r_values_legend.pdf}
	\end{minipage}  
	
	\vskip 0.3cm
	
	\begin{minipage}{0.8\textwidth}
		\centering
		\includegraphics[width=\textwidth]{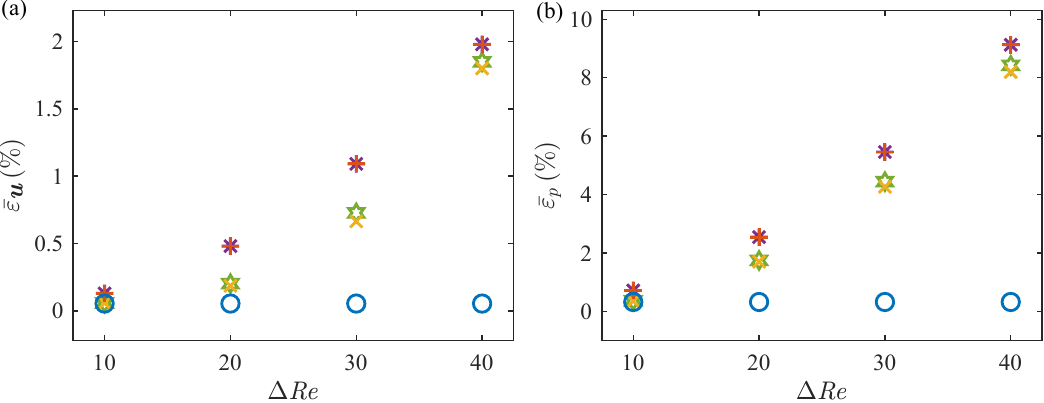}
	\end{minipage} 
	
	\caption{RLEs with respect to system parameter interval $\Delta Re$ for (a) velocity and (b) pressure. (2-point: $Re \in \{$125,135$\}$, $\{$120,140$\}$, $\{$115,145$\}$, $\{$110,150$\}$; 4-point: $Re \in \{$115,125,135,145$\}$, $\{$100,120,140,160$\}$, $\{$85,115,145,175$\}$, $\{$70,110,150,190$\}$).}
	\label{fig:deltaRe}
\end{figure}

\begin{figure}
	\centering
	\begin{minipage}{0.34\textwidth}
		\centering
		\includegraphics[width=\textwidth]{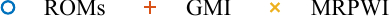}
	\end{minipage}  
	
	\vskip 0.3cm
	
	\begin{minipage}{0.8\textwidth}
		\centering
		\includegraphics[width=\textwidth]{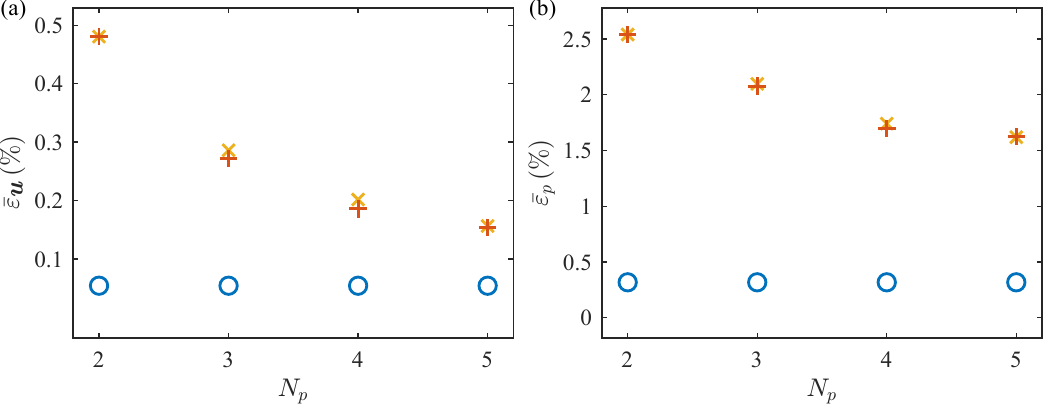}
	\end{minipage} 
	
	\caption{RLEs with respect to number of neighbors $N_p$ for (a) velocity and (b) pressure. ($Re \in \{$120,140$\}$,$\{$100,120,140$\}$, $\{$100,120,140,160$\}$, $\{$80,100,120,140,160$\}$).}
	\label{fig:np}
\end{figure}

POD-Galerkin PROMs targeting $\overline{Re} = 130$ are constructed via GMI and MRPWI, with variations in three key factors: the number of modes $N_r$, the system parameter interval $\Delta Re$, and the number of neighbors $N_p$. Velocity and pressure RLEs decrease monotonically with increasing $N_r$, eventually approaching a plateau, as shown in Fig.~\ref{fig:r}. This trend is observed for both the 2-point and 4-point interpolation methods ($N_p \in \{2,4\}$) at a fixed $\Delta Re = 20$, with $N_r$ ranging from 7 to 21 (step size 2). An increase in the system parameter interval $\Delta Re$ leads to progressively larger RLEs, as depicted in Fig.~\ref{fig:deltaRe}. This trend is consistent for both the 2-point and 4-point interpolation methods ($N_p \in \{2,4\}$) at a fixed $N_r = 13$, with $\Delta Re$ spanning $\{10,20,30,40\}$. RLEs decrease as the number of neighbors $N_p$ increases (Fig.~\ref{fig:np}), indicating improved accuracy with higher-order interpolation. This effect is observed at fixed $N_r = 13$ and $\Delta Re = 20$, with $N_p \in \{2,3,4,5\}$.

These observations indicate that predictive performance of PROMs is enhanced by increasing the number of modes $N_r$ and neighbors $N_p$, while larger system parameter intervals $\Delta Re$ tend to degrade accuracy. Across all tested configurations, PROMs constructed with GMI and MRPWI exhibit similar fidelity. 

\clearpage

\section{Conclusions}\label{sec:conclusions}
POD-Galerkin PROMs, constructed by integrating the POD-Galerkin framework with interpolation methods, are evaluated on the flow over a cylinder benchmark. Specifically, PROMs developed with GMI and the proposed MRPWI achieve high accuracy relative to DNS and standard POD-Galerkin ROMs across a range of configurations, including the number of modes $N_r$, the number of neighbors $N_p$, and the system parameter interval $\Delta Re$. Their accuracy is enhanced with increasing $N_r$ and $N_p$, as well as with decreasing $\Delta Re$. In addition, GMI and the proposed MRPWI exhibit comparable accuracy, while the latter offers improved computational efficiency.

\printcredits

\bibliographystyle{cas-model2-names}

\bibliography{cas-refs}



\end{document}